\documentclass[12pt]{article}
\usepackage{amsfonts}
\usepackage{amssymb,amsmath}
\usepackage{graphicx}
\usepackage{multirow}
\usepackage{hyperref}
\usepackage{float}

\newtheorem{theorem}{Theorem}[section]
\newtheorem{definition}{Definition}[section]
\newtheorem{lemma}{Lemma}[section]

\newtheorem{remark}{Remark}[section]
\newtheorem{condition}{Condition}[section]
\newtheorem{algorithm}{Algorithm}[section]
\allowdisplaybreaks

\begin{document}
\title{A modified subgradient extragradient method for solving the variational inequality problem}
\author{ {\sc Qiao-Li Dong$^a${\thanks{Corresponding author. email: dongql@lsec.cc.ac.cn}},\,\, Aviv Gibali$^b$,\,
Dan Jiang$^a$}\\
\small $^a$Tianjin Key Laboratory for Advanced Signal Processing and College of Science, \\
\small Civil Aviation University of China, Tianjin 300300, China,\\
\small $^b$Department of Mathematics, ORT Braude College,\\
\small  2161002 Karmiel, Israel.\\
}

\date{}

\date{}
\maketitle{} {\bf Abstract.}
The subgradient extragradient method for solving the variational inequality (VI) problem, which is introduced by Censor et al. \cite{CGR}, replaces the second projection onto the feasible set of the VI, in the extragradient method, with a subgradient projection onto some constructible half-space. Since the method has been introduced, many authors proposed extensions and modifications with applications to various problems.

In this paper, we introduce a modified subgradient extragradient method by improving the stepsize of its second step. Convergence of the proposed method is proved under standard and mild conditions and primary numerical experiments illustrate the performance and advantage of this new subgradient extragradient variant.
\\[6pt]
 \\[6pt]

\section{ Introduction}
\vskip 3mm

In this manuscript we are concerned with the variational inequality (VI) problem
of finding a point $x^{\ast }\in \mathcal{H}$ such that%
\begin{equation}
\langle F(x^{\ast }),x-x^{\ast }\rangle \geq 0,\ \text{for all }x\in C,
\label{a}
\end{equation}%
where $C\subseteq \mathcal{H}$ is nonempty, closed and convex set in a real Hilbert
space $\mathcal{H}$, $\langle \cdot ,\cdot \rangle $ denotes the inner product in $\mathcal{H}$,
and $F:\mathcal{H}\rightarrow \mathcal{H}$ is a given mapping.

The VI is a fundamental problem in optimization theory and captures various applications, such as
partial differential equations, optimal control, and mathematical programming (see, for example \cite{Baiocchi-Capelo,FP,Zeidler}).
A vast literature on iterative methods for solving VIs has been published, see for example, \cite{Denisov,Dong1,Dong2,FP,Malitsky,Malitsky-emenov,Noor,YangQZ,Yao,Zhou}. Two special classes of iterative methods which are often used to approximate solutions of the VI problem are presented next.
\vskip 2mm

The first class of methods is the one--step method, also known as projection method, and its iterative step is as follows.
\begin{equation}
x^{k+1}=P_{C}(x^{k}-\alpha _{k}F(x^{k})),
\end{equation}
which is the natural extension of the projected gradient method for optimization problems, originally proposed by Goldstein \cite{Goldstein64} and Levitin and Polyak \cite{lp66}. Under the assumption that $F$ is $\eta-$strongly monotone and $L-$Lipschitz continuous and $\alpha_k\in(0,\frac{2\eta}{L^2})$, the projection method converges to a solution of the VI.
But, if we relax the strong monotonicity assumption to just monotonicity, the situation becomes more complicated, and we may get a divergent sequence independently of the choice of the stepsize $\alpha_k$. To see it, a typical example consists of choosing $C=\mathbb{R}^2$ and $F$ to be a rotation in $\pi/2$, which is certainly monotone and $L$--Lipschitz continuous. The unique solution of the VI \eqref{a} in this case is the origin, but $\{x^k\}_{k=0}^\infty$ gives rise to a sequence satisfying $\|x^{k+1}\|>\|x^k\|$ for all $k\geq0$.
\vskip 2mm

The second class of methods for solving the VI problem is two--steps method. In this class we consider the \textit{extragradient method} introduced by Korpelevich \cite{Korpelevich} and Antipin \cite{Antipin}, which is one of most popular two--steps method, and its iterative step is as follows.
\begin{equation}
\left\{
\begin{array}{l}
y^{k}=P_{C}(x^{k}-\alpha _{k}F(x^{k}))\medskip \\
x^{k+1}=P_{C}(x^{k}-\alpha _{k}F(y^{k}))%
\end{array}%
\right.  \label{EG-1}
\end{equation}%
where $\alpha _{k}\in (0,1/L)$, and $L$ is the Lipschitz constant of $F$, or
$\alpha _{k}$ is updated by the following adaptive procedure
\begin{equation}
\alpha _{k}\Vert F(x^{k})-F(y^{k})\Vert \leq \mu \Vert x^{k}-y^{k}\Vert
,\quad \mu \in (0,1).  \label{a2}
\end{equation}%

The extragradient method has received a great deal of attention and many authors modified and improved it in various ways, see for example \cite{Tseng,Popov}. Here, we focus on one specific extension of He \cite{He} and Sun \cite{Sun}, called the projection and contraction method.
\begin{algorithm}$\left. {}\right. $\textbf{(The projection and contraction method)}
\label{PC}
\begin{equation}
\left\{
\begin{array}{l}
y^{k}=P_{C}(x^{k}-\alpha _{k}F(x^{k}))\medskip \\
x^{k+1}=P_C(x^k-\gamma\rho_k\alpha_k F(y^k))%
\end{array}%
\right.  \label{PC-1}
\end{equation}%
where $\gamma\in (0,2),$  $\alpha_k\in (0,{1}/{L})$ or $\{\alpha_k\}_{k=0}^\infty$ is selected self--adaptively, and
\begin{equation}
 \label{as}
\rho_k:=\frac{\|x^k-y^k\|^2-\alpha_k\langle x^k-y^k, (F(x^k)-F(y^k))\rangle}{\|(x^k-y^k)-\alpha_k (F(x^k)-F(y^k))\|^2},
\end{equation}
\end{algorithm}
\vskip 2mm

The choice of the stepsize is very important since the efficiency of the iterative methods depends heavily on it. Observe that while in the classical extragradient method the stepsize $\alpha_k$ is the same in both projections, here, in the projection and contraction method (\ref{PC-1}), two different stepsizes are used. The numerical experiments presented in \cite{Cai} illustrate that the computational load of the projection and contraction method is about half of that of the extragradient method.
\vskip 2mm

Another observation concerning the extragradient method, is the need to calculate twice
the orthogonal projection onto $C$ per each iteration. So, in case that the set $C$ is not "simple" to project onto it, a minimal distance problem has to be solved (twice) in order to obtain the next iterate, a fact that might affect the efficiency and applicability of the method. As a first step to overcome this obstacle, Censor et al in \cite{CGR}
introduced the \textit{subgradient extragradient method} in which the second projection (\ref{EG-1}) onto $C$ is replaced by a specific subgradient projection which can be easily calculated.
\begin{algorithm}$\left. {}\right. $\textbf{(The subgradient extragradient method)}
\label{SEM}
\begin{equation}
\left\{
\begin{array}{l}
y^{k}=P_{C}(x^{k}-\alpha_k F(x^{k}))\medskip \\
x^{k+1}=P_{T_{k}}(x^{k}-\alpha_k F(y^{k}))%
\end{array}%
\right.
\end{equation}%
where $T_{k}$ is the set defined as%
\begin{equation}
T_{k}:=\{w\in \mathcal{H}\mid \left\langle \left( x^{k}-\alpha_k F(x^{k})\right)
-y^{k},w-y^{k}\right\rangle \leq 0\},  \label{eq:T_K}
\end{equation}
and $\alpha_k\in(0,1/L)$ or $\{\alpha_k\}_{k=0}^\infty$ is selected self--adaptively, that is
$\alpha _{k}=\sigma \rho ^{m_{k}}$, $\alpha >0,$ $\rho \in (0,1)$ and $%
m_{k}$ is the smallest nonnegative integer such that
\begin{equation}
\alpha _{k}\Vert F(x^{k})-F(y^{k})\Vert \leq \mu \Vert x^{k}-y^{k}\Vert
,\quad \mu \in (0,1).  \label{sae}
\end{equation}%
\end{algorithm}

\begin{figure}[H]
\begin{center}
\includegraphics[width=8cm,height=6cm]{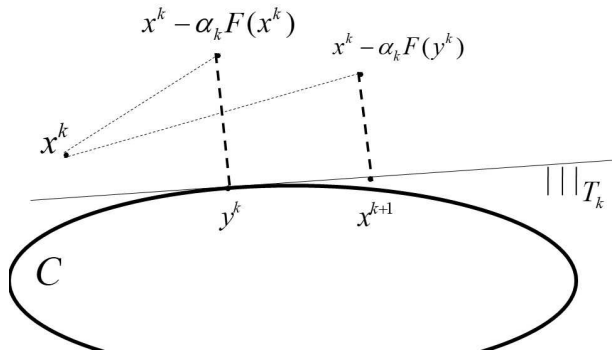}
\end{center}
 \caption{$x^{k+1}$ is a subgradient projection of the point $x^k -\alpha_kF (y^k)$ onto the set $T_k$ }
\label{fig1}
\end{figure}

Figure \ref{fig1} illustrates the iterative step of this algorithm.

Censor et al in \cite{CGR} presented a weak convergence theorem of Algorithm \ref{SEM} with fixed stepsize $\alpha_k=\alpha\in(0,1/L)$, but this result can be easily generalized by using some adaptive step rule as the following theorem shows.

\begin{theorem}
\label{th31}
Given a monotone and $L$-Lipschitz continuous mapping $F:\mathcal{H}\rightarrow\mathcal{H}$. Assume that the solution set of the VI problem (\ref{a}) is nonempty. Then any sequence $\{x^k\}_{k=0}^\infty$ generated by  Algorithm \ref{SEM} satisfies
\begin{equation}
\label{lem32-4}
\|x^{k+1}-x^*\|^2\leq\|x^k-x^*\|^2-(1-\mu^2)\|x^k-y^k\|^2
\end{equation}
and moreover, converges weakly to a solution of the variational inequality problem $(\ref{a})$.
\end{theorem}

Since the inception of the subgradient extragradient method, many authors have proposed various modifications, see for example the results \cite{He-Wu,CGR2,Hieu1}. Kraikaew and Saejung \cite{Kraikaew-Saejung} proposed a Halpern-type variant in order to obtain strong convergence, see also \cite{CGR1}. The subgradient extragradient method is also applied for other probelms than VIs such as multi--valued variational inequality \cite{Fang-Chen}, equilibrium problems \cite{AA,Dang,Dang1} and the split feasibility and fixed point problems \cite{Vinh-Hoai}.
\vskip2mm

So, as mentioned above, the stepsize used extragradient and the subgradient extragradient methods has an essential role in the convergence rate of the two--steps methods, hence it is natural to ask the following question:
\vskip 2mm

\textbf{Is it possible to modify the stepsize in the second step of the subgradient extragradient method in the spirit of He \cite{He} and Sun \cite{Sun}?}
\vskip 2mm

In this paper we provide an affirmative answer to this question by relying on the works of \cite{He,Sun} and introduce a modified subgradient extragradient method which improves the stepsize in the second step of the subgradient extragradient method. To the best of our knowledge, we are not aware of an improvement in the literature. The convergence of the proposed method is proved under standard assumptions and numerical experiment validates its applicability.
\vskip 2mm

The paper is organized as follows: We first recall some basic definitions and results
in Section \ref{sec:Preliminaries}. The modified subgradient extragradient method is presented and analyzed in Section \ref{Sec:MSEM}. Later, in Section \ref{Sec:Num}, some numerical experiments are presented in order to illustrate and compare the performances of the method with other variants. A concluding remarks are given in Section \ref{Sec:Final}.

\section{Preliminaries}\label{sec:Preliminaries}
Let $\mathcal{H}$ be a real Hilbert space with inner product $\langle \cdot
,\cdot \rangle $ and the induced norm $\Vert \cdot \Vert $, and let $D$ be a
nonempty, closed and convex subset of $\mathcal{H}$. We write $%
x^{k}\rightharpoonup x$ to indicate that the sequence $\left\{ x^{k}\right\}
_{k=0}^{\infty }$ converges weakly to $x$. Given a sequence $\left\{ x^{k}\right\}
_{k=0}^{\infty }$, denote by $\omega _{w}(x^{k})$ its weak $\omega $-limit
set, that is, any $x\in \omega _{w}(x^{k})$ such that there exsists a subsequence $%
\left\{ x^{k_{j}}\right\} _{j=0}^{\infty }$of $\left\{ x^{k}\right\}
_{k=0}^{\infty }$ which converges weakly to $x$.
\vskip 2mm

For each point $x\in \mathcal{H},$\ there exists a unique nearest point in $%
D $, denoted by $P_{D}(x)$. That is,%
\begin{equation}
\left\Vert x-P_{D}\left( x\right) \right\Vert \leq \left\Vert x-y\right\Vert
\text{ for all }y\in D.
\end{equation}%
The mapping $P_{D}:\mathcal{H}\rightarrow D$ is called the metric projection
of $\mathcal{H}$ onto $D$. It is well known that $P_{D}$ is a \textit{%
nonexpansive mapping} of $\mathcal{H}$ onto $D$, and further more \textit{%
firmly nonexpansive mapping}. This is captured in the next lemma.

\begin{lemma}
\label{lem21} For any $x,y\in \mathcal{H}$ and $z\in D$, it holds

\begin{itemize}
\item $\Vert P_{D}(x)-P_{D}(y)\Vert\leq \Vert x-y\Vert ;\medskip $

\item $\Vert P_{D}(x)-z\Vert ^{2}\leq \Vert x-z\Vert ^{2}-\Vert
P_{D}(x)-x\Vert ^{2}$;
\end{itemize}
\end{lemma}

The characterization of the metric projection $P_{D}$ \cite[Section 3]%
{Goebel+Reich}, is given in the next lemma.

\begin{lemma}
\label{lem22} Let $x\in \mathcal{H}$ and $z\in D$. Then $z=P_{D}\left(
x\right) $ if and only if
\begin{equation}
P_{D}(x)\in D
\end{equation}%
and%
\begin{equation}
\left\langle x-P_{D}\left( x\right) ,P_{D}\left( x\right) -y\right\rangle
\geq 0\text{ for all }x\in \mathcal{H},\text{ }y\in D.
\end{equation}
\end{lemma}

Given $x \in \mathcal{H}$ and $v \in \mathcal{H}$, $v \neq 0$ and let $T = \{z  \in \mathcal{H}: \langle v, z - x\rangle \leq 0\}.$ Then, for all
$u \in \mathcal{H}$, the projection $P_T (u)$ is defined by
\begin{equation}
\label{projection}
P_T(u)=u-\max\left\{0,\frac{\langle v, u - x\rangle}{\|v\|^2}\right\}v,
\end{equation}
which gives us an explicit formula to find the projection of any point onto a half-space (see \cite{HeSN} for details).

\begin{definition}
The \texttt{normal cone} of $D$ at $v\in D$, denote by $N_{D}\left( v\right)
$ is defined as%
\begin{equation}
N_{D}\left( v\right) :=\{d\in \mathcal{H}\mid \left\langle
d,y-v\right\rangle \leq 0\text{ for all }y\in D\}.  \label{eq:normal_cone}
\end{equation}
\end{definition}

\begin{definition}
Let $B:\mathcal{H}\rightrightarrows 2^{\mathcal{H}}\mathcal{\ }$be a
point-to-set operator defined on a real Hilbert space $\mathcal{H}$. The
operator $B$ is called a \texttt{maximal monotone operator} if $B$ is
\texttt{monotone}, i.e.,%
\begin{equation}
\left\langle u-v,x-y\right\rangle \geq 0\text{ for all }u\in B(x)\text{ and }%
v\in B(y),
\end{equation}%
and the graph $G(B)$ of $B,$%
\begin{equation}
G(B):=\left\{ \left( x,u\right) \in \mathcal{H}\times \mathcal{H}\mid u\in
B(x)\right\} ,
\end{equation}%
is not properly contained in the graph of any other monotone operator.
\end{definition}

It is clear (\cite[Theorem 3]{Rockafellar}) that a monotone mapping $B$ is
maximal if and only if, for any $\left( x,u\right) \in\mathcal{H}\times%
\mathcal{H},$ if $\left\langle u-v,x-y\right\rangle \geq0$ for all $\left(
v,y\right) \in G(B)$, then it follows that $u\in B(x).$

\begin{lemma}
\cite{BC} \label{lem24} Let $D$ be a nonempty, closed and convex subset of
a Hilbert space $\mathcal{H}$. Let $\{x^{k}\}_{k=0}^{\infty }$ be a bounded
sequence which satisfies the following properties:

\begin{itemize}
\item every limit point of $\{x^{k}\}_{k=0}^{\infty }$ lies in $D$;

\item $\lim_{n\rightarrow \infty }\Vert x^{k}-x\Vert $ exists for every $%
x\in D$.
\end{itemize}
Then $\{x^{k}\}_{k=0}^{\infty }$ weakly converges to a point in $D$.
\end{lemma}

\section{The Modified Subgradient Extragradient Method}\label{Sec:MSEM}

In this section, we give a precise statement of our modified subgradient
extragradient method and discuss some of its elementary properties.

\begin{algorithm}
\label{MSEM}$\left. {}\right. $\textbf{(The modified subgradient extragradient method)}
Take $\sigma> 0$, $\rho\in (0, 1)$ and
$\mu\in (0, 1)$.

\textbf{Step 0}: Select a starting point $x^{0}\in\mathcal{H}$
and set $k=0$.

\textbf{Step 1}: Given the current iterate $x^{k},$ compute%
\begin{equation}
y^{k}=P_{C}(x^{k}-\alpha _{k}F(x^{k})),  \label{sStep1}
\end{equation}%
where $\alpha _{k}=\sigma \rho ^{m_{k}}$, $\sigma >0,$ $\rho \in (0,1)$ and $%
m_{k}$ is the smallest nonnegative integer such that
\begin{equation}
\alpha _{k}\Vert F(x^{k})-F(y^{k})\Vert \leq \mu \Vert x^{k}-y^{k}\Vert
,\quad \mu \in (0,1).  \label{sae}
\end{equation}%
If $x^k=y^k$, stop. Otherwise, go to \textbf{Step 2.}

\textbf{Step 2}:
Construct the set
\begin{equation}
T_{k}:=\{w\in \mathcal{H}|\langle (x^{k}-\alpha
_{k}F(x^{k}))-y^{k},w-y^{k}\rangle \leq 0\},  \label{g66}
\end{equation}%
and calculate
\begin{equation}
x^{k+1}=P_{T_{k}}(x^{k}-\gamma\rho_{k}\alpha_kF(y^{k})),  \label{sStep2}
\end{equation}
where $\gamma\in(0,2)$ and
\begin{equation}
\label{step-size}
\aligned
\rho_k:=\frac{\langle x^k-y^k, d(x^k,y^k)\rangle}{\|d(x^k,y^k)\|^2},
\endaligned
\end{equation}
where
\begin{equation}
\label{PC1-OP-d}
d(x^k,y^k)=(x^k-y^k)-\alpha_k(F(x^k)-F(y^k)).
\end{equation}
Set $k\leftarrow (k+1)$ and return to \textbf{Step 1}.
\end{algorithm}
\vskip 2mm

\begin{figure}[H]
\begin{center}
\includegraphics[width=8cm,height=6cm]{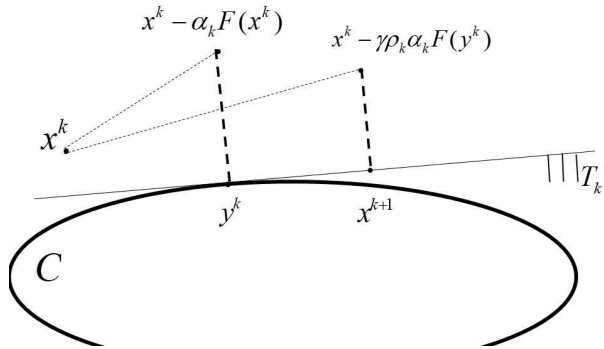}
\end{center}
 \caption{$x^{k+1}$ is a subgradient projection of the point $x^k -\gamma\rho_k\alpha_kF (y^k)$ onto the set $T_k$ }
\label{fig2}
\end{figure}

Figure \ref{fig2} illustrates the iterative step of this algorithm.

Recall that $x^{k}-y^{k}=0$ implies that we are at a solution of the variational inequality. In our convergence theory, we will implicitly assume that this does not occur
after finitely many iterations, so that Algorithm \ref{MSEM} generates an infinite sequence
satisfying, in particular, $x^{k}-y^{k}\neq0$ for all $k\in \mathbb{N}$.
\vskip 2mm

\begin{remark}
We make the following observations for Algorithm \ref{MSEM}.
\begin{itemize}
\item [(1)]
It is easy to see by a simple induction argument from Algorithm \ref{MSEM} that $x^k\in \mathcal{H}$ and $y^k\in C$, which is different from those of the extragradient method \eqref{EG-1}, and projection and contraction method (\ref{PC}).

\item [(2)] The calculation of $\rho_k$ does not add the computational load of the method. The values of $F(x^k)$ and $F(y^k)$ have been obtained in the previous calculation.

\item [(3)]  The projection in \eqref{sStep2} of Algorithm \ref{MSEM} is explicitly computed thanks to formula \eqref{projection}. It is easy to see that $C \subseteq T_k$ for all $k\geq 0$. Indeed, from the definition of $y^k$ and Lemma \ref{lem22}, we obtain $\langle x^k-\alpha_k F(x^k)- y^k, x - y^k\rangle\leq 0$ for all $x \in C$. This together with
the definiton of $T_k$ implies that $C\subseteq T_k$  for all  $k\geq 0$.
\end{itemize}
\end{remark}

\subsection{Convergence Analysis}\label{Sec:Conv}

In this section, we show that Algorithm \ref{MSEM} generates a sequence $\{x^k\}$ which converges weakly to a solution of the variational inequality (\ref{a}). In order to establish this result we assume that the following conditions hold:

\begin{condition}
\label{con:Condition 1.1} The solution set of (\ref{a}), denoted by $SOL(C,F)$,
is nonempty.
\end{condition}

\begin{condition}\label{con:Condition 1.2}
The mapping $F$ is monotone on $\mathcal{H}$, i.e.,
\begin{equation}
\langle F(x)-F(y),x-y\rangle \geq 0,\quad \forall x,y\in \mathcal{H},
\end{equation}
\end{condition}

\begin{condition}\label{con:Condition 1.3}
The mapping $F$ is Lipschitz continuous on $\mathcal{H%
}$ with the Lipschitz constant $L>0$, i.e.,
\begin{equation}
\Vert F(x)-F(y)\Vert \leq L\Vert x-y\Vert ,\quad \forall x,y\in \mathcal{H}.
\end{equation}
\end{condition}

We start our analysis by relying on \cite{Khobotov} shoing that Algorithm \ref{MSEM} is well--defined, meaning that the inner loop of the stepsize calculation in \eqref{sae} is always finite, and that the denominator in the definition of $\alpha_k$ is nonzero.

\begin{lemma} {\rm \cite{Khobotov}} \label{lem41}
The line rule $(\ref{sae})$ is well defined. Besides, $\underline{\alpha}\leq\alpha_k\leq \sigma$, where
$\underline{\alpha}=\min\{\sigma,\frac{\mu\rho}{L}\}.$
\end{lemma}

\begin{lemma}
\label{lem31}
Let $\{\rho_k\}_{k=0}^\infty$ be a sequence defined by (\ref{step-size}). Then under Conditions \ref{con:Condition 1.2} and \ref{con:Condition 1.3}, we have\\
\begin{equation}
\label{le11}
\rho_k\geq\frac{1-\mu}{1+\mu^2}.
\end{equation}
\end{lemma}

\textit{ Proof.}
From the Cauchy-Schwarz inequality and (\ref{sae}), it follows
\begin{equation}
\label{le12}
\aligned
\langle x^k-y^k, d(x^k,y^k)\rangle
&=\langle x^k-y^k,(x^k-y^k)-\alpha_k(F(x^k)-F(y^k))\rangle\\
&=\|x^k-y^k\|^2-\alpha_k\langle x^k-y^k,F(x^k)-F(y^k)\rangle\\
&\geq\|x^k-y^k\|^2-\alpha_k\| x^k-y^k\|\|F(x^k)-F(y^k)\|\\
&\geq(1-\mu)\|x^k-y^k\|^2.\\
\endaligned
\end{equation}
Using Condition \ref{con:Condition 1.2} and \eqref{sae}, we obtain
\begin{equation}
\label{le14}
\aligned
\|d(x^k,y^k)\|^2&=\|x^k-y^k-\alpha_k(F(x^k)-F(y^k))\|^2\\
&=\|x^k-y^k\|^2+\alpha_k^2\|F(x^k)-F(y^k)\|^2-2\alpha_k\langle x^k-y^k,F(x^k)-F(y^k)\rangle\\
&\leq(1+\mu^2)\|x^k-y^k\|^2.\\
\endaligned
\end{equation}
Combining (\ref{le12}) and (\ref{le14}), we obtain (\ref{le11}) and the proof is complete. $\Box$
\vskip 3mm

The contraction property of  Algorithm \ref{MSEM} is presented in the following lemma, which play a key role in the proof of the convergence result.
\begin{lemma}
\label{lem32}
Let sequence $\{x^k\}_{k=0}^\infty$   generated by  Algorithm \ref{MSEM} and $x^*\in SOL(C,F)$.   Then,
under Conditions \ref{con:Condition 1.1}--\ref{con:Condition 1.3}, we have
\begin{equation}
\label{lem32-4}
\|x^{k+1}-x^*\|^2
\leq\|x^k-x^*\|^2-\|(x^k-x^{k+1})-\gamma\rho_kd(x^k,y^k)\|^2
-\gamma(2-\gamma)\rho_k^2\|d(x^k,y^k)\|^2.
\end{equation}
\end{lemma}

\textit{ Proof.}
 By the definition of $x^{k+1}$ and Lemma \ref{lem21}, we have
\begin{equation}
\label{lem32-0}
\aligned
\|x^{k+1}-x^*\|^2
&\leq\|x^k-\gamma\rho_k\alpha_kF(y^k)-x^*\|^2-\|x^k
-\gamma\rho_k\alpha_kF(y^k)-x^{k+1}\|^2\\
&=\|x^k-x^*\|^2-\|x^{k+1}-x^k\|^2-2\gamma\rho_k\alpha_k\langle x^{k+1}-x^*, F(y^k)\rangle.\\
\endaligned
\end{equation}
Since $x^*\in SOL(C,F)$ and $F$ is monotone, we have
$$
\langle F(y^k)-F(x^*),y^k-x^*\rangle\geq0,\quad \forall k\geq0,
$$
which with (\ref{a}) implies
$$
\langle F(y^k),y^k-x^*\rangle\geq0,\quad \forall k\geq0.
$$
So,
\begin{equation}
\label{lem32-1}
\langle F(y^k),x^{k+1}-x^*\rangle\geq\langle F(y^k),x^{k+1}-y^k\rangle.
\end{equation}
By the definition of $T_k$ and $x^{k+1}\in T_k$, we have
$$
\langle (x^{k}-\alpha_{k}F(x^{k}))-y^{k},x^{k+1}-y^{k}\rangle\leq0,
$$
which implies
\begin{equation}
\label{lem32-2}
\langle d(x^k,y^k),x^{k+1}-y^{k}\rangle\leq
\alpha_{k}\langle F(y^{k}),x^{k+1}-y^{k}\rangle.
\end{equation}
Using (\ref{lem32-1}) and (\ref{lem32-2}), we get
\begin{equation}
\label{lem32-3}
\aligned
&-2\gamma\rho_k\alpha_k\langle x^{k+1}-x^*, F(y^k)\rangle\\
&\leq-2\gamma\rho_k\langle x^{k+1}-y^k, d(x^k,y^k)\rangle\\
&=-2\gamma\rho_k\langle x^k-y^k, d(x^k,y^k)\rangle+2\gamma\rho_k\langle x^k-x^{k+1}, d(x^k,y^k)\rangle.\\
\endaligned
\end{equation}
To  the two crossed term in the right hand side of the above formula, we have
\begin{equation}
\label{lem32-31}
-2\gamma\rho_k\langle x^k-y^k, d(x^k,y^k)\rangle=-2\gamma\rho_k^2\|d(x^k,y^k)\|^2,
\end{equation}
and
\begin{equation}
\label{lem32-32}
\aligned
2\gamma\rho_k\langle x^k-x^{k+1}, d(x^k,y^k)\rangle&=-\|(x^k-x^{k+1})-\gamma\rho_kd(x^k,y^k)\|^2\\
&\quad+\|x^k-x^{k+1}\|^2+\gamma^2\rho_k^2\|d(x^k,y^k)\|^2.
\endaligned
\end{equation}
Combining  (\ref{lem32-0}) and (\ref{lem32-3})--\eqref{lem32-32}, we obtain \eqref{lem32-4}.
 $\Box$

\begin{remark}
\begin{itemize}
\item[{\rm(a)}]\
Although $x^{k+1}\notin C$ in Algorithm \ref{MSEM}, by using the definition of $T_k$ and $x^{k+1}\in T_k$, we have the inequality
\begin{equation}
\langle (x^{k}-\alpha_{k}F(x^{k}))-y^{k},x^{k+1}-y^{k}\rangle\leq0,
\end{equation}
which plays a key role in the proof of the contraction inequality \eqref{lem32-4}.
So, we get similar contraction property with \cite[Eq. (4.6)]{Cai}.

\item[{\rm(b)}]\ Employing analysis which are similar  to  those for $\gamma$ after the proof of Theorem 4.1 in \cite{Cai}, we get  that the desirable new iterate $x^{k+1}$ is
updated by (\ref{sStep2}) with $\gamma\in [1, 2)$.
\end{itemize}
\end{remark}

We are now in position to prove our main  convergence result.
\begin{theorem}
\label{th31}
Assume that Conditions \ref{con:Condition 1.1}--\ref{con:Condition 1.3}.  Then
the sequence $\{x^k\}_{k=0}^\infty$   generated by  Algorithm \ref{MSEM} converges weakly to
a solution of the variational inequality problem $(\ref{a})$.
\end{theorem}

\textit{ Proof.}
Fix $x^* \in SOL(C,f )$.
From Condition \ref{con:Condition 1.3}, we have
\begin{equation}
\label{lem32-5}
\aligned
\|d(x^k,y^k)\|&\geq\|x^k-y^k\|-\alpha_k\|F(x^k)-F(y^k)\|\\
&\geq(1-\mu)\|x^k-y^k\|,\\
\endaligned
\end{equation}
Combining  (\ref{le11}), (\ref{lem32-4}) and (\ref{lem32-5}), we get
\begin{equation}
\label{lem32-00}
\aligned
\|x^{k+1}-x^*\|^2
&\leq\|x^k-x^*\|^2-\|(x^k-x^{k+1})-\gamma\rho_kd(x^k,y^k)\|^2\\
&\quad-
\frac{\gamma(2-\gamma)(1-\mu)^3}{1+\mu^2}\|x^k-y^k\|^2.
\endaligned
\end{equation}
From (\ref{lem32-00}), we have
$$
\|x^{k+1}-x^*\|\leq\|x^k-x^*\|,
$$
which implies that the sequence $\{\|x^k-x^*\|\}$ is decreasing and
lower bounded by 0 and thus converges to some finite limit. Moreover, $\{x^k\}_{k=0}^\infty$ is Fej$\acute{e}$r-monotone with respect
to $SOL(C,f)$ and thus is bounded.

From (\ref{lem32-00}) and the existence of $\lim_{k\rightarrow\infty}\|x^k-x^*\|^2$, it follows
\begin{equation}
\label{th31-1}
\sum_{k=0}^{\infty}\|x^{k}-y^k\|\leq+\infty
\end{equation}
which implies
\begin{equation}
\label{th31-2}
\lim_{k\rightarrow\infty}\|x^{k}-y^k\|=0.
\end{equation}

Now, we are to show $\omega_w(x^k)\subseteq SOL(C,F).$ Due to the
boundedness of $\{x^k\}_{k=0}^\infty$, it has at least one weak accumulation point. Let $%
\hat x\in \omega_w(x^k)$. Then there exists a subsequence $\{x^{k_i}\}_{i=0}^\infty$ of $%
\{x^k\}_{k=0}^\infty$ which converges weakly to $\hat x$. From (\ref{th31-2}), it follows
that $\{y^{k_i}\}_{i=0}^\infty$ also converges weakly to $\hat x.$
\vskip 2mm

We will show that $\hat{x}$ is a solution of the variational inequality (\ref%
{a}). Let%
\begin{equation}
A(v)=\left\{
\begin{array}{cc}
F(v)+N_{C}\left( v\right) , & v\in C, \\
\emptyset , & v\notin C\text{,}%
\end{array}%
\right.  \label{cx10}
\end{equation}%
where $N_{C}(v)$ is the normal cone of $C$ at $v\in C$. It is known that $A$
is a maximal monotone operator and $A^{-1}(0)=SOL(C,F)$. If $(v,w)\in G(A)$,
then we have $w-F(v)\in N_{C}(v)$ since $w\in A(v)=F(v)+N_{C}(v)$. Thus it
follows that
\begin{equation}
\langle w-F(v),v-y\rangle \geq 0,\quad y\in C.  \label{g9}
\end{equation}%
Since $y^{k_{i}}\in C$, we have
\begin{equation}
\langle w-F(v),v-y^{k_{i}}\rangle \geq 0.  \label{g10}
\end{equation}

On the other hand, by the definition of $y^k$ and Lemma \ref{lem21}, it
follows that
\begin{equation}  \label{g11}
\langle x^k-\alpha_k F(x^k)-y^k, y^{k}-v\rangle\geq0,
\end{equation}
and consequently,
\begin{equation}  \label{g12}
\left\langle \frac{y^k-x^k}{\alpha_k}+F(x^k),v-y^k\right\rangle\geq0.
\end{equation}
Hence we have
\begin{equation}  \label{g13}
\aligned
&\langle w,v-y^{k_i}\rangle\\
&\geq\langle F(v),v-y^{k_i}\rangle \\
&\geq\langle F(v),v-y^{k_i}\rangle-\Big\langle \frac{y^{k_i}-x^{k_i}}{%
\alpha_{{k_i}}}+F(x^{k_i}),v-y^{k_i}\Big\rangle  \\
&=\langle F(v)-F(y^{k_i}),v-y^{k_i}\rangle+\langle
F(y^{k_i})-F(x^{k_i}),v-y^{k_i}\rangle-\Big\langle \frac{y^{k_i}-x^{k_i}}{%
\alpha_{k_i}},v-y^{k_i}\Big\rangle   \\
&\geq\langle F(y^{k_i})-F(x^{k_i}),v-y^{k_i}\rangle-\Big\langle \frac{%
y^{k_i}-x^{k_i}}{\alpha_{k_i}},v-y^{k_i}\Big\rangle \\
\endaligned
\end{equation}
which implies
\begin{equation}  \label{g14}
\langle w,v-y^{k_i}\rangle\geq\langle F(y^{k_i})-F(x^{k_i}),v-y^{k_i}\rangle-%
\Big\langle \frac{y^{k_i}-x^{k_i}}{\alpha_{k_i}},v-y^{k_i}\Big\rangle. \\
\end{equation}
Taking the limit as $i\rightarrow\infty$ in the above inequality and using
Lemma \ref{lem31}, we obtain
\begin{equation}  \label{g15}
\langle w,v-\hat x\rangle\geq0.
\end{equation}
Since $A$ is a maximal monotone operator, it follows that $\hat x\in
A^{-1}(0) = SOL(C,F)$. So, $\omega_w(x^k)\subseteq SOL(C,F).$

Since $\lim_{k\rightarrow\infty}\|x^{k}-x^*\|$ exists and $%
\omega_w(x^k)\subseteq SOL(C,F)$, using Lemma \ref{lem24}, we conclude that $%
\{x^k\}_{k=0}^\infty$ weakly converges a solution of the variational inequality (\ref{a}%
). This completes the proof. $\Box$

\begin{remark}
\begin{itemize}
\item[{\rm (1)}]\,
With the aid of the proof of Theorem 5.1 in \cite{CGR},
Theorem \ref{th31} will also hold for the fixed stepsize $\alpha_k=\alpha\in(0,1/L).$
\item[{\rm (2)}]\,
The modified subgradient extragradient method can be  generalized  to solve the multi--valued variational inequality in \cite{Fang-Chen} and  the split feasibility and fixed point problems \cite{Vinh-Hoai} since these problems are equivalent with  a inequality problem  or could be easily transformed into an inequality problem. However, the modified subgradient extragradient method couldn't be used directly to solve the equilibrium problems \cite{AA,Dang,Dang1}, which needs further research.
\end{itemize}
\end{remark}

\section{Numerical experiments}\label{Sec:Num}

In this section, we present a numerical example to compare the modified subgradient extragradient method (Algorithm \ref{MSEM}) with the subgradient extragradient method (Algorithm \ref{SEM}) and the projection and contraction method (Algorithm \ref{PC}).
\vskip 0.2cm

Consider the linear operator $Ax: = M x + q$, which is taken from \cite{Harker} and has been considered by many authors for numerical experiments, see, for example \cite{Hieu, Solodov}, where
\begin{equation}
M= BB^T+ S + D,
\end{equation}
and $B$ is an $m\times m$ matrix, $S$ is an $m\times m$ skew-symmetric matrix, $D$ is an $m\times m$
diagonal matrix, whose diagonal entries are nonnegative (so $M$ is positive semidefinite),
$q$ is a vector in $\mathbb{R}^m$. The feasible set $C\subset \mathbb{R}^m$  is closed and convex and defined as
\begin{equation}
C: = \{x\in \mathbb{R}^m \mid Qx\leq b\},
\end{equation}
where $Q$ is an $l\times m$ matrix and $b$ is a nonnegative vector. It is clear that $A$ is monotone and $L$--Lipschitz continuous with $L =\|M\|$ (hence uniformly continuous). For $ q= 0$, the solution set $SOL(C,A)= \{0\}$.

Just as in \cite{Hieu}, we randomly choose the starting points $x^1\in[0,1]^m$ in  Algorithms \ref{PC}, \ref{SEM} and  \ref{MSEM}. We choose the stopping criterion as $\|x^k\|\leq \epsilon = 0.005$ and the parameters $\sigma=7.55,\, \rho=0.5,\, \mu=0.85$ and $\gamma=1.99$. The size $l = 100$ and $m = 5, 10, 20, 30, 40, 50, 60, 70, 80$. The matrices $B,S,D$ and the vector $b$ are generated randomly.
\vskip 2mm

In Table \ref{table1}, we denote by ``Iter." the number of iterations and ``InIt." the number of total iterations of finding suitable $\alpha_k$ in (\ref{sae}).
 \vskip 0.2cm

\begin{table}[H]\caption{Comparison of Algorithms \ref{PC}, \ref{SEM} and  \ref{MSEM}}\label{table1}
{\footnotesize
\begin{center}
\begin{tabular}{l|c c c c c c c c c c c c }
\hline
&&\multicolumn{1}{c}{Iter. }&&&\multicolumn{2}{c}{InIt.}&&&\multicolumn{2}{c}{CPU in second}&\\
\cline{2-4} \cline{6-8}\cline{10-12}
 $m$ & Alg. \ref{PC} & Alg. \ref{SEM} &Alg. \ref{MSEM} &&Alg. \ref{PC} & Alg. \ref{SEM} &Alg. \ref{MSEM} &&Alg. \ref{PC} & Alg. \ref{SEM} &Alg. \ref{MSEM}\\
 \hline
 \hline
 5   &24&84&21&&166&487&146&&0.8594&0.2500&0.3438\\
 10  &59&149&60&&502&1022&512&&0.6250&0.1719&0.1719\\
 20  &99&1145&199&&962&10290&2139&&32.3750&2.9844&2\\
 30  &484&1137&485&&5714&10692&5727&&2.1094&1.1563&1.2188\\
 40  &733&2814&648&&9004&28063&8281&&98.2188&24.0781&4.9531\\
 50  &1234&4809&1526&&16218&51843&20606&&643.7344&8.5313&2.6563\\
 60  &1431&7475&712&&19276&82188&9968&&257.0781&28.2344&8.5156\\
 70  &--&13016&2350&&--&155821&35167&&--&45.1563&5.6406\\
 80  &2894&12270&2200&&41915&145878&32425&&187.6563&26.2500&8.2969\\
\hline
\hline
\end{tabular}
\end{center}
}
\end{table}

Table \ref{table1}, shows that Algorithm \ref{MSEM} highly improves Algorithm \ref{SEM} with respect to the number of iterations and CPU time. ``Iter." and ``InIt." are almost the same for Algorithm \ref{MSEM} and \ref{PC}, however,  Algorithm \ref{MSEM} needs less CPU time than  Algorithm \ref{PC} since the projection onto $C$ is more complicated than projection onto $T_k$.

\section{Final Remarks}\label{Sec:Final}

In this article, we propose a modified subgradient extragradient method for solving the VI problem by improving
the stepsize in the second step of the subgradient extragradient method. Under standard assumptions weak convergence of the proposed method is presented. Preliminary numerical experiments indicate that our method does greatly outperform the subgradient extragradient method. Since there are other two-steps variants for solving the VI problem, such as \cite{Popov,Malitsky-emenov1}, in which only one evaluation of the mapping $F$ is needed per each iteration, it is thus natural to apply our stepsize strategy to this and other two--steps methods, and this is one of our future research topics.
\vskip 4mm

\noindent
{\bf Acknowledgements.}
\vskip 2mm

{The authors express their thanks to the two anonymous referees, whose careful readings and suggestions led to improvements in the presentation of the results.

The first author is supported by National Natural Science Foundation of China (No. 71602144)  and Open Fund of Tianjin Key Lab for Advanced Signal Processing (No. 2016ASP-TJ01). The third author is supported by the EU FP7 IRSES program STREVCOMS, grant no. PIRSES-GA-2013-612669.
}
\vskip 4mm

 \end{document}